\documentclass[12pt]{article}
\usepackage{amssymb, amsmath}

\def\R{\mathbb{R}}
\def\C{\mathbb{C}}
\def\N{\mathbb{N}}

\def\I{1\!\!1}

\newtheorem{theorem}{Theorem}
\newtheorem{corollary}[theorem]{Corollary}
\newtheorem{lemma}[theorem]{Lemma}
\newtheorem{proposition}[theorem]{Proposition}
\newtheorem{definition}[theorem]{Definition}

\newtheorem{remark}[theorem]{Remark}

\begin{document}

\author{\textbf{Maria Jo\~{a}o Oliveira} \\
{\small Universidade Aberta, P 1269-001 Lisbon, Portugal}\\
{\small CMAF, University of Lisbon, P 1649-003 Lisbon, Portugal}\\
{\small oliveira@cii.fc.ul.pt} \and
\textbf{Jos\'e Lu{\'\i}s da Silva} \\
{\small DME, University of Madeira, P 9000-390 Funchal, Portugal}\\
{\small CCM, University of Madeira, P 9000-390 Funchal, Portugal}\\
{\small luis@uma.pt} \and
\textbf{Ludwig Streit} \\
{\small Forschungszentrum BiBoS, Universit{\"a}t Bielefeld, D 33501
Bielefeld,
Germany}\\
{\small CCM, University of Madeira, P 9000-390 Funchal, Portugal}\\
{\small streit@physik.uni-bielefeld.de}}

\title{Intersection local times of independent fractional Brownian motions 
as generalized white noise functionals}

\date{}
\maketitle

\vspace*{-0.5cm}

\begin{abstract}
In this work we present expansions of intersection local times of fractional 
Brownian motions in $\R^d$, for any dimension $d\geq 1$, with arbitrary Hurst 
coefficients in $(0,1)^d$. The expansions are in terms of Wick powers of white 
noises (corresponding to multiple Wiener integrals), being well-defined in the 
sense of generalized white noise functionals. As an application of our 
approach, a sufficient condition on $d$ for the existence of intersection 
local times in $L^2$ is derived, extending the results in \cite{NL07} to 
different and more general Hurst coefficients.
\end{abstract}


\medskip

\noindent
\textbf{2000 AMS Classification:} 60H40, 60G15, 60J55, 28C20, 46F25, 82D60

\section{Introduction}

In the recent years the fractional Brownian motion has become an object of 
intense study, namely, due to its special properties, such as short/long range 
dependence and self-similarity, yielding its proper and natural uses in 
several applications in different fields (e.g.~mathematical finances 
\cite{MO08}, telecommunications engineering \cite{NRT03}). 

Besides its own specific properties, the intersection properties of fractional 
Brownian motion paths have been studied by many authors as well, see e.g.~the 
works done by Gradinaru et al.~\cite{GRV03}, Nualart et al.~\cite{NH07}, 
\cite{NH05}, Rosen \cite{R87}, and the references therein.

One may consider intersections of sample paths with themselves, as in 
\cite{DOS08} and references therein, or with other independent fractional 
Brownian motions, as in \cite{NL07}. 

This work concerns the latter standpoint. Within the white noise analysis 
framework (Section \ref{Section2}), a first purpose of this work is an 
extension of the results presented in \cite{AOS01} to two $d$-dimensional 
independent fractional Brownian motions $\mathbf{B}_{H_1}$ and 
$\mathbf{B}_{H_2}$ with different Hurst coefficients, $H_1$ and $H_2$. 
Technically, this approach has the advantage that the underlying probability 
space does not depend on any Hurst coefficient under consideration. As a 
consequence, one may analyze the intersection local time of any two 
independent fractional Brownian motions, without any restriction on the 
corresponding Hurst coefficients. 

From the viewpoint of applications to physics, this absence of restrictions on 
the Hurst coefficients under consideration is meaningful to widen the 
modelling of polymers towards polymers molecules handling different types of 
polymers.

For low dimensions, that is, either for $d=1$ or for $d=2$, the white noise 
analysis framework allows the definition of the intersection local time of any 
two independent fractional Brownian motions $\mathbf{B}_{H_i}$ in terms of an 
integral over a Donsker's $\delta$-function
$$
L\equiv \int d^2t\, \delta(\mathbf{B}_{H_1}(t_1) - \mathbf{B}_{H_2}(t_2)),
$$
intended to sum up the contributions from each pair of moments of time $t_1$, 
$t_2$ for which the fractional Brownian motions $\mathbf{B}_{H_i}$ arrive at 
the same point. 

A rigorous definition, such as, e.g., through a sequence of Gaussians 
approximating the $\delta$-function,
\[
(2\pi\varepsilon)^{-d/2}\exp\left(-\frac{|x|^{2}}{2\varepsilon}\right),\quad
\varepsilon>0,
\]
will make $L$ increasingly singular, and various ``renormalizations'' 
have to be done as the dimension $d$ increases. Of course, besides the 
dimension of the space, the type of ``renormalizations'' needed depends as 
well on the Hurst coefficients $H_i\in\left(0,1\right)^d$ being considered. For 
$d>2$ with $1/{\max_jH_{1,j}}+1/{\max_jH_{2,j}}\leq d$, the expectation diverges 
in the limit and must be subtracted. Depending on the values of 
$\max_jH_{i,j}$, further kernel terms must be also subtracted (Theorem 
\ref{Th9}).   

In this work we are particularly interested in the chaos decomposition of $L$. 
We expand $L$ in terms of Wick powers \cite{HKPS93} of white noise,
an expansion which corresponds to that in terms of multiple Wiener integrals
when one considers the Wiener process as the fundamental random variable.
This allows us to derive the kernels for $L$. Due to the local structure of the 
Wick powers, the kernel functions are relatively simple and exhibit clearly 
the dimension dependence singularities of $L$ (Proposition \ref{Proposition2}).
For comparison, we also calculate the regularized kernel functions
corresponding to the Gaussian $\delta$-sequence mentioned above (Theorem 
\ref{Th9}).

As an application of this approach, in Theorem \ref{Th11} we derive a 
sufficient condition for the existence of the intersection local times in 
$L^2$, extending the results obtained in \cite{NL07} to different and more 
general Hurst coefficients. 

\section{Gaussian white noise calculus}\label{Section2}

In this section we briefly recall the concepts and results of white
noise analysis used throughout this work (for a detailed explanation see 
e.g.~\cite{BeKo88}, \cite{HKPS93}, \cite{HOUZ96}, \cite{Kuo96}, \cite{Ob94}). 

\subsection{Fractional Brownian motion}

The starting point of white noise analysis for the construction of two 
independent $d$-dimensional, $d\geq 1$, fractional Brownian motions is the 
real Gelfand triple 
\[
S_{2d}(\R)\subset L_{2d}^2(\R)\subset S_{2d}'(\R),
\]
where $L_{2d}^{2}(\R):=L^2(\R,\R^{2d})$ is the real Hilbert space of all 
vector valued square integrable functions with respect to the Lebesgue 
measure on $\R$, and $S_{2d}(\R)$, $S_{2d}'(\R)$ are the Schwartz spaces of the 
vector valued test functions and tempered distributions, respectively. We shall
denote the $L^2_{2d}(\R)$-norm by $|\cdot|_{2d}$ (or if there is no risk of 
confusion simply by $|\cdot|$) and the dual pairing between $S_{2d}'(\R)$ and 
$S_{2d}(\R)$ by $\left\langle \cdot,\cdot\right\rangle_{2d}$, or simply by 
$\left\langle \cdot,\cdot\right\rangle$, which is defined as the bilinear 
extension of the inner product on $L_{2d}^2(\R)$, i.e.,
\[
\langle\mathbf{g},\mathbf{f}\rangle_{2d}= \sum_{i=1}^{2d} \int_{\R} dx\,g_i(x)f_i(x), 
\]
for all $\mathbf{g}=(g_1,..., g_{2d})\in L^2_{2d}(\R)$ and all
$\mathbf{f}=(f_1,..., f_{2d})\in S_{2d}(\R)$. By the Minlos theorem, there is
a unique probability measure $\mu$ on the $\sigma$-algebra $\mathcal{B}$ 
generated by the cylinder sets on $S'_{2d}(\R)$ with characteristic function 
given by
\[
C(\mathbf{f}):=\int_{S'_{2d}(\R)}d\mu(\vec\omega)\,e^{i\left\langle\vec\omega,\mathbf{f}\right\rangle }
=e^{-\frac{1}{2}|\mathbf{f}|^{2}},\quad\mathbf{f}\in S_{2d}(\R).
\]
In this way we have defined the white noise measure space 
$(S_{2d}'(\R),\mathcal{B},\mu)$. 

To construct two independent $d$-dimensional fractional Brownian motions we 
shall consider a $2d$-tuple of independent Gaussian white noises
\[
\vec\omega:=(\vec\omega_1,\vec\omega_2),\quad
\vec\omega_i=(\omega_{i,1},...,\omega_{i,d}), i=1,2.
\]
Within this formalism, a version of a $d$-dimensional Wiener Brownian motion 
is given by 
\[
\mathbf{B}(t):=\left(\langle\omega_1,\I_{[0,t]}\rangle,...,
\langle\omega_d,\I_{[0,t]}\rangle\right),\quad (\omega_1,...,\omega_d)\in S_d'(\R),
\]
where $\I_A$ denotes the indicator function of a set $A$ and 
$\langle\cdot,\cdot\rangle=\langle\cdot,\cdot\rangle_1$. For an arbitrary 
$d$-dimensional Hurst parameter $H=(H_1,...,H_d)\in\left(0,1\right)^d$, a 
version of a $d$-dimensional fractional Brownian motion is given by
\[
\mathbf{B}_H(t):=\left(\langle\omega_1,M_{H_1}\I_{[0,t]}\rangle,...,
\langle\omega_d,M_{H_d}\I_{[0,t]}\rangle\right),\quad 
(\omega_1,...,\omega_d)\in S_d'(\R),
\]
where, for a 1-dimensional Hurst parameter $H\in\left(0,1\right)$ and for a 
generic real valued function $f$, 
\begin{eqnarray*}
(M_Hf)(x):=\left\{ 
\begin{array}{cl}
&\!\!\!\!\displaystyle\frac{(\frac12-H)K_H}{\Gamma\left(H+\frac12\right)}\lim_{\varepsilon\to 0^+}\int_\varepsilon^{\infty}dy\,\frac{f(x)-f(x+y)}{y^{\frac32-H}},\ 
H\in\left(0,1/2\right)\ (*)\\
& \\
&\!\!\!\!f(x),\ H=\frac12\\
&  \\
&\!\!\!\!\displaystyle\frac{K_H}{\Gamma\left(H-\frac12\right)}
\int_x^{\infty}dy\,f(y)(y-x)^{H-\frac32},\ H\in\left(1/2,1\right)\ (**)\
\end{array}
\right. ,
\end{eqnarray*}
provided the limit in ($*$) exists for almost all $x\in\R$ and the 
integral in ($**$) exists for all $x\in\R$ (for more details see 
e.g.~\cite{B03} and \cite{PT00} and the references therein). Independently of 
the case under consideration, the normalizing constant $K_H$ is given by
\[
K_H=\Gamma\left(H+\frac12\right)
\left(\frac{1}{2H}+\int_0^{\infty}ds\,\left((1+s)^{H-\frac12}-s^{H-\frac12}\right)\right)^{-\frac12}.
\]
There are several examples of functions $f$ for which $M_Hf$
exists for any $H\in(0,1)$, namely, for $f=\I_{[0,t]}$ with $t>0$ or for 
$f\in S_1(\R)$. For more details and proofs see e.g.~\cite{B03}, \cite{BHOZ08}, \cite{M08}, \cite{PT03}, and the references therein.

\subsection{Hida distributions and characterization results}

Let us now consider the complex Hilbert space 
$(L^2):=L^{2}(S'_{2d}(\R),\mathcal{B},\mu)$. For simplicity one introduces the 
notation
\[
\mathbf{n} =(n_1, \cdots ,n_d)\in \N^d,\quad n = \sum_{i = 1}^d n_i,\quad 
\mathbf{n}! = \prod_{i = 1}^d n_i!.
\]
The space $(L^{2})$ is canonically isomorphic to the symmetric Fock space of 
symmetric square integrable functions,
$$
(L^{2})\simeq 
\Big(\bigoplus_{k = 0}^\infty \mathrm{Sym}\, L^2(\R^k, k!d^kx)\Big)^{\otimes 2d},
$$
which leads to the chaos expansion of the elements in $(L^{2})$, 
\begin{eqnarray*} 
F(\vec\omega_1,\vec\omega_2) &=& \sum_{\mathbf{m}}\sum_{\mathbf{k}}\langle:\vec\omega_1^{\otimes\mathbf{m}}: \otimes :\vec\omega_2^{\otimes\mathbf{k}}:,\mathbf{f}_{\mathbf{m},\mathbf{k}}\rangle\\
&=&\sum_{\mathbf{m}}\sum_{\mathbf{k}}
\left\langle\bigotimes_{i=1}^d:\omega_{1,i}^{\otimes m_i}:\otimes\bigotimes_{j=1}^d:\omega_{2,j}^{\otimes m_j}:,\mathbf{f}_{\mathbf{m},\mathbf{k}}\right\rangle,
\end{eqnarray*}
with kernel functions $\mathbf{f}_{\mathbf{m},\mathbf{k}}$ in the Fock space, 
that is, square integrable functions of the $m+k$ arguments and symmetric in 
each $m_i$-, $k_j$-tuple. 

To proceed further we have to consider a Gelfand triple around the space 
$(L^{2})$. We will use the space $(S)^*$ of Hida distributions (or generalized 
Brownian functionals) and the corresponding Gelfand triple 
$(S)\subset (L^{2})\subset(S)^*$. Here $(S)$ is the space of white noise test 
functions such that its dual space (with respect to $(L^{2})$) is the space 
$(S)^*$. Instead of reproducing the explicit construction of 
$\left(S\right)^*$ (see e.g.~\cite{HKPS93}), in Theorem \ref{08Prop1} below we 
characterize this space through its $S$-transform. We recall that given a 
$\mathbf{f}=(\mathbf{f}_1,\mathbf{f}_2)\in S_{2d}(\R)$, and the Wick exponential
\[
:\exp(\langle\vec\omega,\mathbf{f}\rangle):\,:=
\sum_{\mathbf{m}}\sum_{\mathbf{k}}\frac{1}{\mathbf{m}!\mathbf{k}!}
\langle:\vec\omega_1^{\otimes\mathbf{m}}:\otimes:\vec\omega_2^{\otimes\mathbf{k}}:,
\mathbf{f}_1^{\otimes\mathbf{m}}\otimes\mathbf{f}_2^{\otimes\mathbf{k}}\rangle
=C(\mathbf{f})e^{\langle\vec\omega,\mathbf{f}\rangle_{2d}},
\]
we define the $S$-transform of a $\Phi\in \left(S\right)^*$ by
\begin{equation}
S\Phi(\mathbf{f}):= 
\left\langle\!\left\langle \Phi,:\exp(\left\langle \cdot,\mathbf{f}
\right\rangle):\right\rangle\!\right\rangle,\quad \forall\,\mathbf{f}\in S_{2d}(\R).\label{08eq6}
\end{equation}
Here $\left\langle\!\left\langle\cdot ,\cdot\right\rangle\!\right\rangle$ 
denotes the dual pairing between $\left(S\right)^*$ and 
$\left(S\right)$ which is defined as the bilinear extension of 
the sesquilinear inner product on $(L^2)$. We observe that the multilinear 
expansion of (\ref{08eq6}),
\[
S\Phi(\mathbf{f}):= \sum_{\mathbf{m}}\sum_{\mathbf{k}}\langle F_{\mathbf{m},\mathbf{k}},\mathbf{f}_1^{\otimes\mathbf{m}}\otimes\mathbf{f}_2^{\otimes\mathbf{k}}\rangle,
\] 
extends the chaos expansion to $\Phi\in \left(S\right)^*$ with distribution valued kernels $F_{\mathbf{m},\mathbf{k}}$ such that
\begin{equation}
\left\langle\!\left\langle\Phi,\varphi\right\rangle\!\right\rangle
=\sum_{\mathbf{m}}\sum_{\mathbf{k}}\mathbf{m}!\mathbf{k}!\langle F_{\mathbf{m},\mathbf{k}},\varphi_{\mathbf{m},\mathbf{k}}\rangle,\label{norm}
\end{equation}
for every generalized test function $\varphi\in(S)$ with kernel functions 
$\varphi_{\mathbf{m},\mathbf{k}}$.

In order to characterize the space $\left(S\right)^*$ through its 
$S$-transform we need the following definition.

\begin{definition}
\label{Def1}A function $F:S_{2d}(\R)\rightarrow \C$ is called a 
$U$-functional whenever\newline
1. for every $\mathbf{f}_1,\mathbf{f}_2\in S_{2d}(\R)$ the mapping 
$\R\ni \lambda \longmapsto F(\lambda \mathbf{f}_1+\mathbf{f}_2)$ has an 
entire extension to $\lambda \in \C$,\newline
2. there are constants $K_1,K_2>0$ such that 
\[
\left| F(z\mathbf{f})\right| \leq K_1e^{K_2\left| z\right| ^2\left\|
\mathbf{f}\right\| ^2} ,\quad \forall \,z\in \C,\mathbf{f}\in S_{2d}(\R)
\]
for some continuous norm $\left\| \cdot \right\|$ on $S_{2d}(\R)$.
\end{definition}

We are now ready to state the aforementioned characterization result.

\begin{theorem}
\label{08Prop1}{\rm (\cite{KLPSW96}, \cite{PS91})} The $S$-transform defines a 
bijection between the space $\left(S\right)^*$ and the space of 
$U$-functionals.
\end{theorem}

As a consequence of Theorem \ref{08Prop1} one may derive the next two 
statements. The first one concerns the convergence of sequences of 
Hida distributions and the second one the Bochner integration of families of 
distributions of the same type (for more details and proofs see 
e.g.~\cite{HKPS93}, \cite{KLPSW96}, \cite{PS91}).

\begin{corollary}
\label{Corollary2}Let $\left( \Phi _n\right)_{n\in\N}$ be a sequence in 
$\left(S\right)^*$ such that
\begin{description}     
\item[{\it (i)}] for all $\mathbf{f}\in S_{2d}(\R)$, 
$\left((S\Phi_n)(\mathbf{f})\right)_{n\in\N}$ is a Cauchy sequence in $\C$, 
\item[{\it (ii)}] there are constants $K_1,K_2>0$ such that for some 
continuous norm $\left\| \cdot \right\|$ on $S_{2d}(\R)$ one has  
\[
\left| (S\Phi_n)(z\mathbf{f})\right| \leq K_1e^{K_2\left| z\right| ^2\left\|
\mathbf{f}\right\| ^2} ,\quad \forall \,z\in \C,\mathbf{f}\in S_{2d}(\R),n\in\N.
\]
\end{description}
Then $\left(\Phi _n\right)_{n\in\N}$ converges strongly in 
$\left(S\right)^*$ to a unique Hida distribution. 
\end{corollary}

\begin{corollary}
\label{Corollary1} Let $(\Omega, \mathcal{B}, m)$ be a measure space and 
$\lambda\mapsto \Phi_\lambda$ be a mapping from $\Omega$ to $(S)^*$. We assume that the $S$-transform of $\Phi_\lambda$ fulfills the following two properties:
\begin{description}     
\item[{\it (i)}] the mapping 
$\lambda\mapsto (S\Phi_\lambda)(\mathbf{f})$ is measurable for every $\mathbf{f}\in S_{2d}(\R)$,
\item[{\it (ii)}] the $S\Phi_\lambda$ obeys a $U$-estimate
\[
|(S\Phi_\lambda)(z\mathbf{f})| \leq C_1(\lambda) e^{C_2(\lambda) |z|^2 
\Vert\mathbf{f}\Vert^2},\quad z\in\C,\mathbf{f}\in S_{2d}(\R),
\]
for some continuous norm $\Vert\cdot\Vert$ on $S_{2d}(\R)$ and for some
$C_1\in L^1(\Omega,m)$, $C_2\in L^\infty(\Omega,m)$.
\end{description}
Then 
\[
\int_\Omega dm(\lambda)\, \Phi_\lambda\in (S)^*
\]
and
\[
S\left(\int_\Omega dm(\lambda)\,\Phi_\lambda\right) (\mathbf{f}) =
\int_\Omega dm(\lambda)\,\left(S\Phi_\lambda\right)(\mathbf{f}).
\] 
\end{corollary}

\section{Chaos expansions}

Let us now consider two independent $d$-dimensional fractional Brownian motions 
$\mathbf{B}_{H_1}(t)$ and $\mathbf{B}_{H_2}(t)$ with Hurst multiparameters 
$H_1=(H_{1,1},...,H_{1,d})$ and $H_2=(H_{2,1},...,H_{2,d})$, respectively. That 
is, given a $2d$-tuple of independent white noises 
$(\omega_{1,1},...,\omega_{1,d},\omega_{2,1},...,\omega_{2,d})$,
\[
\mathbf{B}_{H_i}(t):=\left(\langle\omega_{i,1},M_{H_{i,1}}\I_{[0,t]}\rangle,...,
\langle\omega_{i,d},M_{H_{i,d}}\I_{[0,t]}\rangle\right),\quad i=1,2.
\] 

\begin{proposition} For each $t$ and $s$ strictly positive real numbers the 
Bochner integral 
\[
\delta (\mathbf{B}_{H_1}(t)-\mathbf{B}_{H_2}(s)):=\left(\frac 1{2\pi }\right)^d
\int_{\R^d}d\mathbf{\lambda}\,e^{i\mathbf{\lambda}(\mathbf{B}_{H_1}(t)-\mathbf{B}_{H_2}(s))}
\]
is a Hida distribution with $S$-transform given by
\begin{eqnarray}
&&S\delta (\mathbf{B}_{H_1}(t)-\mathbf{B}_{H_2}(s))(\mathbf{f})\nonumber\\
&=&\left(\frac{1}{\sqrt{2\pi}}\right)^d
\prod_{j=1}^d\frac{1}{\sqrt{t^{2H_{1,j}}+s^{2H_{2,j}}}}\cdot\label{eq4}\\
&&\cdot e^{-\frac12\sum_{j=1}^d\frac{1}{t^{2H_{1,j}}+s^{2H_{2,j}}}
\left(\int_{\R}dx\,\left(f_{1,j}(x)(M_{H_{1,j}}\I_{\left[0,t\right]})(x)-f_{2,j}(x)(M_{H_{2,j}}\I_{\left[0,s\right]})(x)\right)\right)^2},\nonumber
\end{eqnarray}
for all $\mathbf{f}=(f_{1,1},...,f_{1,d},f_{2,1},...,f_{2,d})\in S_{2d}(\R)$.
\end{proposition}

\noindent
\textbf{Proof.} The proof of this result follows from an application of 
Corollary \ref{Corollary1} to the $S$-transform of the integrand function 
\[
\Phi(\vec\omega_1, \vec\omega_2):= 
e^{i\mathbf{\lambda}(\mathbf{B}_{H_1}(t)-\mathbf{B}_{H_2}(s))},\quad 
\vec\omega_i=(\omega_{i,1},...,\omega_{i,d}), i=1,2,
\]
with respect to the Lebesgue measure on $\R^d$. For this purpose we begin by 
observing that since the fractional Brownian motions are independent one has
\[
S\Phi(\mathbf{f})=Se^{i\mathbf{\lambda}\mathbf{B}_{H_1}(t)}(\mathbf{f}_1)\cdot
Se^{-i\mathbf{\lambda}\mathbf{B}_{H_2}(s)}(\mathbf{f}_2)
\]
for every $\mathbf{f}=(\mathbf{f}_1,\mathbf{f}_2)\in S_{2d}(\R)$, 
$\mathbf{f}_1:=(f_{1,1},...,f_{1,d})$, $\mathbf{f}_2:=(f_{2,1},...,f_{2,d})$. 
Hence, according e.g.~to \cite{HKPS93}, for all 
$\mathbf{\lambda}=(\lambda_1,...,\lambda_d)\in \R^d$ we obtain
\begin{equation}
\!\!\!S\Phi(\mathbf{f})
=\prod_{j=1}^de^{i\lambda_j\int_{\R}dx\,
\left(f_{1,j}(x)(M_{H_{1,j}}\I_{\left[0,t\right]})(x)-f_{2,j}(x)(M_{H_{2,j}}\I_{\left[0,s\right]})(x)\right)}e^{-\frac12 \lambda_j^2(t^{2H_{1,j}}+s^{2H_{2,j}})}\label{eq2},
\end{equation}
which clearly fulfills the measurability condition. Moreover, for all $z\in\C$ 
we find 
\begin{eqnarray*}
&&\left|S\Phi(z\mathbf{f})\right|\\
&=&\prod_{j=1}^de^{-\frac14 \lambda_j^2(t^{2H_{1,j}}+s^{2H_{2,j}})}\cdot \\
&&\cdot\prod_{j=1}^d\left|e^{-\frac14 \lambda_j^2(t^{2H_{1,j}}+s^{2H_{2,j}})+iz\lambda_j\int_{\R}dx\,
\left(f_{1,j}(x)(M_{H_{1,j}}\I_{\left[0,t\right]})(x)-f_{2,j}(x)(M_{H_{2,j}}\I_{\left[0,s\right]})(x)\right)}\right|\\
&\leq&\prod_{j=1}^de^{-\frac14 \lambda_j^2(t^{2H_{1,j}}+s^{2H_{2,j}})}\cdot\\
&&\cdot\prod_{j=1}^de^{-\frac14 \lambda_j^2(t^{2H_{1,j}}+s^{2H_{2,j}})+|z||\lambda_j|
\left|\int_{\R}dx\,
\left(f_{1,j}(x)(M_{H_{1,j}}\I_{\left[0,t\right]})(x)-f_{2,j}(x)(M_{H_{2,j}}\I_{\left[0,s\right]})(x)\right)\right|},
\end{eqnarray*}
where, for each $j=1,...,d$, the corresponding term in the second product is 
bounded by
\[
\exp\left(\frac{|z|^2}{t^{2H_{1,j}}+s^{2H_{2,j}}}\left(\int_{\R}dx\,\left(f_{1,j}(x)(M_{H_{1,j}}\I_{\left[0,t\right]})(x)-f_{2,j}(x)(M_{H_{2,j}}\I_{\left[0,s\right]})(x)\right)\right)^2\right),
\]
because
\begin{eqnarray*}
&&-\frac14 \lambda_j^2(t^{2H_{1,j}}+s^{2H_{2,j}})\\
&&+|z||\lambda_j|
\left|\int_{\R}dx\,
\left(f_{1,j}(x)(M_{H_{1,j}}\I_{\left[0,t\right]})(x)-f_{2,j}(x)(M_{H_{2,j}}\I_{\left[0,s\right]})(x)\right)\right|\\
&=&-\left(\frac{|z|}{\sqrt{t^{2H_{1,j}}+s^{2H_{2,j}}}}\left|\int_{\R}dx\,
\left(f_{1,j}(x)(M_{H_{1,j}}\I_{\left[0,t\right]})(x)-f_{2,j}(x)(M_{H_{2,j}}\I_{\left[0,s\right]})(x)\right)\right|\right.\\
&&\left.-\frac{|\lambda_j|}{2}\sqrt{t^{2H_{1,j}}+s^{2H_{2,j}}}\right)^2\\
&&+\frac{|z|^2}{t^{2H_{1,j}}+s^{2H_{2,j}}}\left(\int_{\R}dx\,\left(f_{1,j}(x)(M_{H_{1,j}}\I_{\left[0,t\right]})(x)-f_{2,j}(x)(M_{H_{2,j}}\I_{\left[0,s\right]})(x)\right)\right)^2.
\end{eqnarray*}
As a result,
\begin{eqnarray*}
\left|S\Phi(z\mathbf{f})\right|
&\leq&e^{-\frac14 \sum_{j=1}^d\lambda_j^2(t^{2H_{1,j}}+s^{2H_{2,j}})}\cdot\\
&&\cdot e^{|z|^2\sum_{j=1}^d\frac{1}{t^{2H_{1,j}}+s^{2H_{2,j}}}\left(\int_{\R}dx\,\left(f_{1,j}(x)(M_{H_{1,j}}\I_{\left[0,t\right]})(x)-f_{2,j}(x)(M_{H_{2,j}}\I_{\left[0,s\right]})(x)\right)\right)^2},
\end{eqnarray*}
where, as a function of $\mathbf{\lambda}$, the first exponential is 
integrable on $\R^d$ and the second exponential is constant. 

An application of the result mentioned above completes the proof. In 
particular, it yields (\ref{eq4}) by integrating (\ref{eq2}) over 
$\mathbf{\lambda}$.\hfill$\blacksquare \medskip$

In order to proceed further the next result shows to be very useful. It
improves the estimate obtained in \cite[Theorem 2.3]{B03} towards the 
characterization results stated in Corollaries \ref{Corollary2} and 
\ref{Corollary1}.

\begin{lemma} (\cite{DOS08})
\label{Lemma} Let $H\in (0,1)$ and $f\in S_1(\R)$ be given. 
There is a non-negative constant $C_H$ independent of $f$ such that
\[
\left|\int_{\R}dx\,f(x)(M_H\I_{\left[0,t\right]})(x)\right|
\leq C_H t\left(\sup_{x\in\R}|f(x)|+ \sup_{x\in\R}|f'(x)|
+|f|\right)
\]
for all $t>0$.
\end{lemma}

In particular, the use of Lemma \ref{Lemma} allows to state the next result on 
intersection local times $L_{H_1,H_2}$ as well as on their subtracted 
counterparts $L_{H_1,H_2}^{(N)}$. There, and throughout the rest of this work 
as well, given a $H=(H_1,...,H_d)\in\left(0,1\right)^d$ we shall use the 
notation
\[
\bar H:= \max_{j=1,...,d}H_j.
\]

\begin{theorem}
\label{Theorem1} Let $T>0$ be given. For any pair of integer numbers $d\geq 1$, 
$N\geq 0$ and for any pair of Hurst multiparameters $H_1,H_2\in (0,1)^d$ 
such that 
\[
\max\{\bar H_1,\bar H_2\}\left(N+\frac{d}{2}-\frac{1}{2\min\{\bar H_1,\bar H_2\}}\right)<N+\frac12, 
\]
the Bochner integral
\[
L_{H_1,H_2}^{(N)}:=\int_0^T dt\int_0^Tds\,\delta^{(N)}(\mathbf{B}_{H_1}(t)-\mathbf{B}_{H_2}(s))
\]
is a Hida distribution.
\end{theorem}

\noindent
\textbf{Proof.} To prove this result we shall again use Corollary 
\ref{Corollary1} with respect to the Lebesgue measure on $\left[0,T\right]^2$. 
For this purpose let us denote the truncated exponential series by
\[
\exp_N (x) := \sum_{n = N}^\infty {{x^n}\over {n!}}.
\] 
It follows from (\ref{eq4}) that for every $t,s>0$ the $S$-transform of 
$\delta^{(N)}(\mathbf{B}_{H_1}(t)-\mathbf{B}_{H_2}(s))$ is given by 
\begin{eqnarray}
&&S\delta^{(N)}(\mathbf{B}_{H_1}(t)-\mathbf{B}_{H_2}(s))(\mathbf{f})\label{eq6}\\
&=&\left(\frac{1}{\sqrt{2\pi}}\right)^d
\prod_{j=1}^d\frac{1}{\sqrt{t^{2H_{1,j}}+s^{2H_{2,j}}}}
\exp_N\left(-\frac12\sum_{j=1}^d\frac{1}{t^{2H_{1,j}}+s^{2H_{2,j}}}\right.\cdot\nonumber\\
&&\cdot\left.\left(\int_{\R}dx\,\left(f_{1,j}(x)(M_{H_{1,j}}\I_{\left[0,t\right]})(x)-f_{2,j}(x)(M_{H_{2,j}}\I_{\left[0,s\right]})(x)\right)\right)^2\right),\nonumber
\end{eqnarray}
which is a measurable function.

In order to check the boundedness condition, on $S_{2d}(\R)$ let us consider 
the norm $\Vert\cdot\Vert$ defined for all 
$\mathbf{f}=(f_1,...,f_{2d})\in S_{2d}(\R)$ by
\begin{equation}
\Vert\mathbf{f}\Vert:= \left(\sum_{i=1}^{2d}\left(\sup_{x\in\R}|f_i(x)|
+ \sup_{x\in\R}|f_i'(x)| + |f_i|\right)^2\right)^{\frac12}.\label{08equ1}
\end{equation}
We observe that on $S_1(\R)$ this norm reduces to the continuous norm
\[
\Vert f\Vert=\sup_{x\in\R}|f(x)|+\sup_{x\in\R}|f'(x)| + |f|,\quad f\in S_1(\R),
\]
which implies the continuity of the norm (\ref{08equ1}) for higher dimensions. 

By Lemma \ref{Lemma}, for each $j=1,...,d$, we obtain
\begin{eqnarray*}
&&\left(\int_{\R}dx\,\left(f_{1,j}(x)(M_{H_{1,j}}\I_{\left[0,t\right]})(x)-f_{2,j}(x)(M_{H_{2,j}}\I_{\left[0,s\right]})(x)\right)\right)^2\\
&\leq&2\left(\int_{\R}dx\,f_{1,j}(x)(M_{H_{1,j}}\I_{\left[0,t\right]})(x)\right)^2+2 \left(\int_{\R}dx\,f_{2,j}(x)(M_{H_{2,j}}\I_{\left[0,s\right]})(x)\right)^2\\
&\leq&2t^2C_{H_{1,j}}^2\Vert f_{1,j}\Vert^2+2s^2C_{H_{2,j}}^2\Vert f_{2,j}\Vert^2,
\end{eqnarray*}
and thus, for all $z\in\C$ and all $\mathbf{f}\in S_{2d}(\R)$,
\begin{eqnarray*}
&&\!\!\!\left|S(\delta^{(N)}(\mathbf{B}_{H_1}(t)-\mathbf{B}_{H_2}(s)))(z\mathbf{f})\right|\\
&\leq&\!\!\!\left(\frac{1}{\sqrt{2\pi}}\right)^d
\prod_{j=1}^d\frac{1}{\sqrt{t^{2H_{1,j}}+s^{2H_{2,j}}}}
\exp_N\left(|z|^2C_{H_1,H_2}^2\frac{t^2+s^2}{t^{2H_{(1)}}+s^{2H_{(2)}}}\Vert\mathbf{f}\Vert^2\right)
\end{eqnarray*}
with $C_{H_1,H_2}:=\max\{C_{H_{1,j}},C_{H_{2,j}}:j=1,...,d\}$ and 
\begin{equation*}
H_{(1)}:=
\begin{cases}
\bar H_1=\displaystyle\max_{j=1,...,d}H_{1,j},&0<t\leq 1\\
\displaystyle\min_{j=1,...,d}H_{1,j},&t>1 
\end{cases},\ 
H_{(2)}:=
\begin{cases}
\bar H_2=\displaystyle\max_{j=1,...,d}H_{2,j},&0<s\leq 1\\
\displaystyle\min_{j=1,...,d}H_{2,j},&s>1
\end{cases}.
\end{equation*}
Therefore, for $0<t,s\leq 1$ one has
\[
\frac{t^2+s^2}{t^{2H_{(1)}}+s^{2H_{(2)}}}=\frac{t^2+s^2}{t^{2\bar H_1}+s^{2\bar H_2}}\leq 1,
\]
and thus
\[
\exp_N\left(|z|^2C_{H_1,H_2}^2\frac{t^2+s^2}{t^{2\bar H_1}+s^{2\bar H_2}}\Vert\mathbf{f}\Vert^2\right)\leq
\left(\frac{t^2+s^2}{t^{2\bar H_1}+s^{2\bar H_2}}\right)^Ne^{|z|^2C_{H_1,H_2}^2\Vert\mathbf{f}\Vert^2};
\]
while either for $t>1$ or for $s>1$ one finds
\begin{eqnarray*}
&&\exp_N\left(|z|^2C_{H_1,H_2}^2\frac{t^2+s^2}{t^{2H_{(1)}}+s^{2H_{(2)}}}\Vert\mathbf{f}\Vert^2\right)\\
&\leq& 
\left(\frac{t^2+s^2}{t^{2H_{(1)}}+s^{2H_{(2)}}}\right)^Ne^{|z|^2C_{H_1,H_2}^2\left(\frac{t^2+s^2}{t^{2H_{(1)}}+s^{2H_{(2)}}}+N\right)\Vert\mathbf{f}\Vert^2}.
\end{eqnarray*}
As a consequence, independently of $T$ being smaller or greater than 1 there 
is always a function $C=C(t,s)>0$ bounded on $\left[0,T\right]^2$ such that
\begin{eqnarray}
&&\left|S(\delta^{(N)}(\mathbf{B}_{H_1}(t)-\mathbf{B}_{H_2}(s)))(z\mathbf{f})\right|\nonumber\\
&\leq& \left(\frac{1}{\sqrt{2\pi}}\right)^d\prod_{j=1}^d\frac{1}{\sqrt{t^{2H_{1,j}}+s^{2H_{2,j}}}}\left(\frac{t^2+s^2}{t^{2H_{(1)}}+s^{2H_{(2)}}}\right)^Ne^{|z|^2C_{H_1,H_2}^2C(t,s)\Vert\mathbf{f}\Vert^2}.\label{eq3}
\end{eqnarray}  
The proof then amounts to prove the integrability on $\left[0,T\right]^2$
of the expression
\[
\prod_{j=1}^d\frac{1}{\sqrt{t^{2H_{1,j}}+s^{2H_{2,j}}}}\left(\frac{t^2+s^2}{t^{2H_{(1)}}+s^{2H_{(2)}}}\right)^N
\]
appearing in (\ref{eq3}).
For this purpose one observes that due to the singular point at the origin 
this expression is integrable on $\left[0,T\right]^2$ if and only if it is 
integrable on $\left[0,1\right]^2$. As shown in the Appendix (Lemma 
\ref{Lemma2}), this occurs whenever 
\[
2\max\{\bar H_1,\bar H_2\}\left(N+\frac{d}{2}-\frac{1}{2\min\{\bar H_1,\bar H_2\}}\right)-2N<1. 
\]
The proof is then completed by an application of Corollary \ref{Corollary1}.\hfill$\blacksquare \medskip$ 

As a consequence, one may derive the chaos expansion for the (truncated) local 
times $L_{H_1,H_2}^{(N)}$.

\begin{proposition}
\label{Proposition2} Under the conditions of Theorem \ref{Theorem1}, 
$L_{H_1,H_2}^{(N)}$ has the chaos expansion 
\[
L_{H_1,H_2}^{(N)}(\vec\omega_1,\vec\omega_2)=\sum_{\mathbf{m}}\sum_{\mathbf{k}}\langle:\vec\omega_1^{\otimes\mathbf{m}}: \otimes :\vec\omega_2^{\otimes\mathbf{k}}:,F_{H_1,H_2,\mathbf{m},\mathbf{k}}\rangle
\]
where the kernel functions $F_{H_1,H_2,\mathbf{m},\mathbf{k}}$ are 
given by 
\begin{eqnarray*}
&&F_{H_1,H_2,\mathbf{m},\mathbf{k}}=\\
&&\left(\frac{1}{\pi}\right)^{\frac{d}{2}}\frac{(-1)^{\frac{m+3k}{2}}}{\left(\frac{\mathbf{m}+\mathbf{k}}{2}\right)!}
\left(\frac12\right)^{\frac{m+k+d}{2}}\binom{\mathbf{m}+\mathbf{k}}{\mathbf{m}}
\int_0^Tdt\int_0^Tds\prod_{j=1}^d\left(\frac{1}{t^{2H_{1,j}}+s^{2H_{2,j}}}\right)^{\frac{m_j+k_j+1}{2}}\cdot\\
&&\cdot\bigotimes_{j=1}^d\left((M_{H_{1,j}}\I_{\left[0,t\right]})^{\otimes m_j}
\otimes(M_{H_{2,j}}\I_{\left[0,s\right]})^{\otimes k_j}\right)
\end{eqnarray*}
for each $\mathbf{m}=(m_1,\ldots,m_d)$ and each $\mathbf{k}=(k_1,\ldots,k_d)$ such that
$m+k\geq 2N$ and all sums $m_j+k_j$, $j=1,...,d$, are even numbers. All other 
kernel functions $F_{H_1,H_2,\mathbf{m},\mathbf{k}}$ are identically 
equal to zero.
\end{proposition} 

\noindent
\textbf{Proof.} According to Corollary \ref{Corollary1}, the $S$-transform of 
the (truncated) local time $L_{H_1,H_2}^{(N)}$ is obtained by integrating 
(\ref{eq6}) over $\left[0,T\right]^2$. Hence, given a 
$\mathbf{f}=(f_{1,1},...,f_{1,d},f_{2,1},...,f_{2,d})\in S_{2d}(\R)$ one has
\begin{eqnarray}
&&SL_{H_1,H_2}^{(N)}(\mathbf{f})\nonumber\\
&=&\left(\frac{1}{\sqrt{2\pi}}\right)^d\int_0^Tdt\int_0^Tds\prod_{j=1}^d\frac{1}{\sqrt{t^{2H_{1,j}}+s^{2H_{2,j}}}}\cdot\nonumber\\
&&\cdot\sum_{n=N}^\infty\frac{(-1)^n}{2^nn!}\sum_{{n_1,\cdots ,n_d}\atop{n_1 +\cdots + n_d = n}}\frac{n!}{n_1!\cdots n_d!}\prod_{j=1}^d\left(\frac{1}{t^{2H_{1,j}}+s^{2H_{2,j}}}\right)^{n_j}\cdot\nonumber\\
&&\cdot\left(\int_{\R}dx\,\left(f_{1,j}(x)(M_{H_{1,j}}\I_{\left[0,t\right]})(x)-f_{2,j}(x)(M_{H_{2,j}}\I_{\left[0,s\right]})(x)\right)\right)^{2n_j}\label{exp}
\end{eqnarray}
with (\ref{exp}) being equal to
\begin{eqnarray*}
&&\sum_{{m_j,k_j}\atop{m_j+k_j=2n_j}}(-1)^{k_j}\binom{2n_j}{m_j}
\left(\int_{\R}dx\,f_{1,j}(x)(M_{H_{1,j}}\I_{\left[0,t\right]})(x)\right)^{m_j}\cdot\\
&&\qquad\qquad\cdot\left(\int_{\R}dx\,f_{2,j}(x)(M_{H_{2,j}}\I_{\left[0,s\right]})(x)\right)^{k_j}.
\end{eqnarray*}
From these calculations follow the equality
\begin{eqnarray*}
SL_{H_1,H_2}^{(N)}(\mathbf{f})
&=&\left(\frac{1}{\pi}\right)^{\frac{d}{2}}\int_0^Tdt\int_0^Tds
\sum_{n = N}^\infty\sum_{{n_1,\cdots ,n_d}\atop{n_1 +\cdots + n_d = n}}\sum_{{m_1,\ldots,m_d, k_1,\ldots,k_d}\atop{m_j + k_j = 2n_j, j = 1, \ldots, d}}\\
&&\left\{\prod_{j=1}^d\frac{(-1)^{\frac{m_j + 3k_j}{2}}}{\left(\frac{m_j+k_j}{2}\right)!}\left(\frac{1}{2(t^{2H_{1,j}}+s^{2H_2,j})}\right)^{\frac{m_j + k_j+1}{2}}\right\}\cdot\\
&&\cdot\left\{\prod_{j=1}^d\binom{m_j+k_j}{m_j}\left(\int_{\R}dx\,f_{1,j}(x)(M_{H_{1,j}}\I_{\left[0,t\right]})(x)\right)^{m_j}\right.\cdot\\
&&\qquad\qquad\qquad\quad\cdot\left.\left(\int_{\R}dx\,f_{2,j}(x)(M_{H_{2,j}}\I_{\left[0,s\right]})(x)\right)^{k_j}\right\},
\end{eqnarray*}
which is equivalent to
\begin{eqnarray*}
&&\left(\frac{1}{\pi}\right)^{\frac{d}{2}}\int_0^Tdt\int_0^Tds
\sum_{{{\mathbf{m}, \mathbf{k}}\atop{m + k\geq 2N}}\atop{m_j + k_j\, even, j=1,...,d}}
\frac{(-1)^{\frac{m+3k}{2}}}{\left(\frac{\mathbf{m}+\mathbf{k}}{2}\right)!}
\left(\frac12\right)^{\frac{m+k+d}{2}}\binom{\mathbf{m}+\mathbf{k}}{\mathbf{m}}\cdot\\
&&\cdot\prod_{j=1}^d\left(\frac{1}{t^{2H_{1,j}}+s^{2H_{2,j}}}\right)^{\frac{m_j+k_j+1}{2}}\left(\int_{\R}dx\,f_{1,j}(x)(M_{H_{1,j}}\I_{\left[0,t\right]})(x)\right)^{m_j}\\
&&\qquad\left(\int_{\R}dx\,f_{2,j}(x)(M_{H_{2,j}}\I_{\left[0,s\right]})(x)\right)^{k_j}.
\end{eqnarray*}
Comparing with the general form of the chaos expansion
\[
\sum_{\mathbf{m}}\sum_{\mathbf{k}}\langle:\vec\omega_1^{\otimes\mathbf{m}}: \otimes :\vec\omega_2^{\otimes\mathbf{k}}:,F_{H_1,H_2,\mathbf{m},\mathbf{k}}\rangle,
\]
one concludes that the kernels $F_{H_1,H_2,\mathbf{m},\mathbf{k}}$ 
vanish whenever either there is a $j=1,...,d$ such that $m_j+k_j$ is an odd 
number or $m+k<2N$, while for all other cases
\begin{eqnarray*}
&&F_{H_1,H_2,\mathbf{m},\mathbf{k}}=\\
&&\left(\frac{1}{\pi}\right)^{\frac{d}{2}}\frac{(-1)^{\frac{m+3k}{2}}}{\left(\frac{\mathbf{m}+\mathbf{k}}{2}\right)!}
\left(\frac12\right)^{\frac{m+k+d}{2}}\binom{\mathbf{m}+\mathbf{k}}{\mathbf{m}}
\int_0^Tdt\int_0^Tds\prod_{j=1}^d\left(\frac{1}{t^{2H_{1,j}}+s^{2H_{2,j}}}\right)^{\frac{m_j+k_j+1}{2}}\cdot\\
&&\cdot\bigotimes_{j=1}^d\left((M_{H_{1,j}}\I_{\left[0,t\right]})^{\otimes m_j}
\otimes(M_{H_{2,j}}\I_{\left[0,s\right]})^{\otimes k_j}\right).
\end{eqnarray*}
\hfill$\blacksquare \medskip$ 

Theorem \ref{Theorem1} shows that for $d=1$ or $d=2$ all intersection local 
times $L_{H_1,H_2}$ are well-defined for all possible Hurst multiparameters 
$H_1$, $H_2$ in $(0,1)^d$. For $d>2$, intersection local times are well-defined 
only for $1/{\bar H_1}+1/{\bar H_2}>d$. Under these conditions, Proposition 
\ref{Proposition2} in addition yields 
\[
\mathbb{E}_\mu(L_{H_1,H_2})=F_{H_1,H_2,0,0}=\left(\frac{1}{\sqrt{2\pi}}\right)^d\int_0^Tdt\int_0^Tds\,\prod_{j=1}^d\frac{1}{\sqrt{t^{2H_{1,j}}+ s^{2H_{2,j}}}}.
\]

Informally speaking, for $1/{\bar H_1}+1/{\bar H_2}\leq d$ with $d>2$, the 
local times only become well-defined once subtracted the divergent terms. This 
``renormalization'' procedure motivates the study of a regularization. As a 
computationally simple regularization we discuss
\[
L_{H_1,H_2,\varepsilon}:=\int_0^Tdt\int_0^Tds\,\delta_\varepsilon(\mathbf{B}_{H_1}(t)-\mathbf{B}_{H_2}(s)),\quad \varepsilon>0,
\]
where
\[
\delta_\varepsilon(\mathbf{B}_{H_1}(t)-\mathbf{B}_{H_2}(s)):=
\left(\frac 1{\sqrt{2\pi\varepsilon}}\right)^de^{-\frac{(\mathbf{B}_{H_1}(t)-\mathbf{B}_{H_2}(s))^2}{2\varepsilon}}.
\]

\begin{theorem}
\label{Th9}
Let $\varepsilon >0$ be given. For all $H_1,H_2\in (0,1)^d$ and all 
dimensions $d\geq 1$ the intersection local time $L_{H_1,H_2,\varepsilon}$ is a 
Hida distribution with kernel functions given by
\begin{eqnarray*}
&&F_{H_1,H_2,\varepsilon, \mathbf{m},\mathbf{k}}=
\left(\frac{1}{\pi}\right)^{\frac{d}{2}}\frac{(-1)^{\frac{m+3k}{2}}}{\left(\frac{\mathbf{m}+\mathbf{k}}{2}\right)!}
\left(\frac12\right)^{\frac{m+k+d}{2}}\binom{\mathbf{m}+\mathbf{k}}{\mathbf{m}}\\
&&\int_0^T\!\!dt\!\!\int_0^T\!\!ds\prod_{j=1}^d\left(\frac{1}{\varepsilon+t^{2H_{1,j}}+s^{2H_{2,j}}}\right)^{\frac{m_j+k_j+1}{2}}\!\!\!\!
\bigotimes_{j=1}^d\left((M_{H_{1,j}}\I_{\left[0,t\right]})^{\otimes m_j}\!
\otimes\!(M_{H_{2,j}}\I_{\left[0,s\right]})^{\otimes k_j}\right)
\end{eqnarray*}
for all $\mathbf{m}=(m_1,...,m_d),\mathbf{k}=(k_1,...,k_d)\in\N_0^d$ such that 
all sums $m_i+k_j$, $j=1,...,d$, are even numbers, and 
$F_{H_1,H_2,\varepsilon, \mathbf{m},\mathbf{k}}\equiv 0$ if at least one of the 
sums $m_i+k_i$ is an odd number. Moreover, if $\max\{\bar H_1,\bar H_2\}\left(N+\frac{d}{2}-\frac{1}{2\min\{\bar H_1,\bar H_2\}}\right)<N+\frac12$, then when 
$\varepsilon$ tends to zero the (truncated) intersection local time 
$L^{(N)}_{H_1,H_2,\varepsilon}$ converges strongly in $(S)^*$ to the (truncated) 
local time $L_{H_1,H_2}^{(N)}$.  
\end{theorem}

\noindent
\textbf{Proof.} As before, the first part of the proof follows from the 
Corollary \ref{Corollary1} with respect to the Lebesgue measure on 
$\left[0,T\right]^2$. By the definition of the $S$-transform, for all 
$\mathbf{f}=(f_{1,1},...,f_{1,d},f_{2,1},...,f_{2,d})\in S_{2d}(\R)$ one finds 
\begin{eqnarray*}
&&S\delta_\varepsilon(\mathbf{B}_{H_1}(t)-\mathbf{B}_{H_2}(s))(\mathbf{f})\\
&=&\prod_{j=1}^d\frac{1}{\sqrt{2\pi(\varepsilon+t^{2H_{1,j}}+s^{2H_{2,j}})}}\cdot\\
&&\cdot e^{-\frac12\sum_{j=1}^d\frac{1}{\varepsilon+t^{2H_{1,j}}+s^{2H_{2,j}}}
\left(\int_{\R}dx\,\left(f_{1,j}(x)(M_{H_{1,j}}\I_{\left[0,t\right]})(x)-f_{2,j}(x)(M_{H_{2,j}}\I_{\left[0,s\right]})(x)\right)\right)^2},
\end{eqnarray*}
which is measurable. Hence, similarly to the proof of Theorem \ref{Theorem1},
an application of Lemma \ref{Lemma} yields for all $z\in\C$ and all 
$\mathbf{f}\in S_{2d}(\R)$ 
\begin{eqnarray*}
&&\!\!\!\left|S(\delta_\varepsilon(\mathbf{B}_{H_1}(t)-\mathbf{B}_{H_2}(s)))(z\mathbf{f})\right|\\
&\leq&\!\!\!\prod_{j=1}^d\frac{1}{\sqrt{2\pi(\varepsilon+t^{2H_{1,j}}+s^{2H_{2,j}})}}e^{|z|^2C_{H_1,H_2}^2\frac{t^2+s^2}{\varepsilon+t^{2H_{(1)}}+s^{2H_{(2)}}}\Vert\mathbf{f}\Vert^2}
\end{eqnarray*}
with $\frac{t^2+s^2}{\varepsilon+t^{2H_{(1)}}+s^{2H_{(2)}}}$ bounded on 
$\left[0,T\right]^2$ and $\prod_{j=1}^d\frac{1}{\sqrt{2\pi(\varepsilon+t^{2H_{1,j}}+s^{2H_{2,j}})}}$ integrable on $\left[0,T\right]^2$. By 
Corollary \ref{Corollary1}, one may then conclude that 
$L_{H_1,H_2,\varepsilon}\in (S)^*$ and, moreover, for every 
$\mathbf{f}=(f_{1,1},...,f_{1,d},f_{2,1},...,f_{2,d})\in S_{2d}(\R)$,
\begin{eqnarray*}
SL_{H_1,H_2,\varepsilon}(\mathbf{f})\!\!&=&\!\!\int_0^Tdt\int_0^Tds\,S\delta_\varepsilon(\mathbf{B}_{H_1}(t)-\mathbf{B}_{H_2}(s))(\mathbf{f})\\
&=&\left(\frac{1}{\pi}\right)^{\frac{d}{2}}\int_0^Tdt\int_0^Tds\!\!\!\!\!
\sum_{{\mathbf{m}, \mathbf{k}}\atop{m_j + k_j\, even, j=1,...,d}}\!\!\!\!\!
\frac{(-1)^{\frac{m+3k}{2}}}{\left(\frac{\mathbf{m}+\mathbf{k}}{2}\right)!}
\left(\frac12\right)^{\frac{m+k+d}{2}}\binom{\mathbf{m}+\mathbf{k}}{\mathbf{m}}\cdot\\
&&\cdot\prod_{j=1}^d\left(\frac{1}{\varepsilon+t^{2H_{1,j}}+s^{2H_{2,j}}}\right)^{\frac{m_j+k_j+1}{2}}\left(\int_{\R}dx\,f_{1,j}(x)(M_{H_{1,j}}\I_{\left[0,t\right]})(x)\right)^{m_j}\\
&&\qquad\left(\int_{\R}dx\,f_{2,j}(x)(M_{H_{2,j}}\I_{\left[0,s\right]})(x)\right)^{k_j}.
\end{eqnarray*}
As in the proof of Proposition \ref{Proposition2}, it follows from the latter 
expression that the kernels $F_{H_1,H_2,\varepsilon, \mathbf{m},\mathbf{k}}$ 
appearing in the chaos expansion of $L_{H_1,H_2,\varepsilon}$ vanish if at least 
one of the $m_i+k_i$ in $\mathbf{m}+\mathbf{k}=(m_1+k_1,...,m_d+k_d)$ is an 
odd number, otherwise they are given by
\begin{eqnarray*}
&&F_{H_1,H_2,\varepsilon, \mathbf{m},\mathbf{k}}=
\left(\frac{1}{\pi}\right)^{\frac{d}{2}}\frac{(-1)^{\frac{m+3k}{2}}}{\left(\frac{\mathbf{m}+\mathbf{k}}{2}\right)!}
\left(\frac12\right)^{\frac{m+k+d}{2}}\binom{\mathbf{m}+\mathbf{k}}{\mathbf{m}}\\
&&\int_0^T\!\!dt\!\!\int_0^T\!\!ds\prod_{j=1}^d\left(\frac{1}{\varepsilon+t^{2H_{1,j}}+s^{2H_{2,j}}}\right)^{\frac{m_j+k_j+1}{2}}\!\!\!\!
\bigotimes_{j=1}^d\left((M_{H_{1,j}}\I_{\left[0,t\right]})^{\otimes m_j}\!
\otimes\!(M_{H_{2,j}}\I_{\left[0,s\right]})^{\otimes k_j}\right).
\end{eqnarray*}
To complete the proof amounts to check the convergence. For this purpose we 
shall use Corollary \ref{Corollary2}. Since
\[
SL^{(N)}_{H_1,H_2,\varepsilon}(\mathbf{f})=\int_0^Tdt\int_0^Tds\,S\delta^{(N)}_\varepsilon(\mathbf{B}_{H_1}(t)-\mathbf{B}_{H_2}(s))(\mathbf{f}),
\]
for every $z\in\C$ and every $\mathbf{f}\in S_{2d}(\R)$, a similar procedure 
used to prove Theorem \ref{Theorem1} yields
\begin{eqnarray*}
\left|SL^{(N)}_{H_1,H_2,\varepsilon}(z\mathbf{f})\right|&\leq&\int_0^Tdt\int_0^Tds\,\left|S\delta^{(N)}_\varepsilon(\mathbf{B}_{H_1}(t)-\mathbf{B}_{H_2}(s))(z\mathbf{f})\right|\\
&\leq&\left(\frac{1}{\sqrt{2\pi}}\right)^d
e^{|z|^2\Vert\mathbf{f}\Vert^2C_{H_1,H_2}^2\sup_{t,s\in\left[0,T\right]}C(t,s)}\cdot\\
&&\cdot\int_0^Tdt\int_0^Tds\,\prod_{j=1}^d\frac{1}{\sqrt{t^{2H_{1,j}}+s^{2H_{2,j}}}}\left(\frac{t^2+s^2}{t^{2H_{(1)}}+s^{2H_{(2)}}}\right)^N,
\end{eqnarray*}
showing the boundedness condition. Furthermore, we have
\begin{eqnarray*}
&&\left|S\delta^{(N)}_\varepsilon(\mathbf{B}_{H_1}(t)-\mathbf{B}_{H_2}(s))(\mathbf{f})\right|\leq\\
&&\left(\frac{1}{\sqrt{2\pi}}\right)^d\prod_{j=1}^d\frac{1}{\sqrt{t^{2H_{1,j}}+s^{2H_{2,j}}}}\left(\frac{t^2+s^2}{t^{2H_{(1)}}+s^{2H_{(2)}}}\right)^Ne^{C_{H_1,H_2}^2\Vert\mathbf{f}\Vert^2\sup_{t,s\in\left[0,T\right]}C(t,s)},
\end{eqnarray*}
which allows the use of the Lebesgue dominated convergence theorem to infer 
the other condition needed for the application of Corollary \ref{Corollary2}.\hfill$\blacksquare \medskip$

Given any pair of Hurst multiparameters $H_1, H_2\in(0,1)^d$, $d\geq 1$, such 
that $d<1/{\bar H_1}+1/{\bar H_2}$, according to the convergence result stated 
in Theorem \ref{Th9}, for any $\mathbf{f}\in S_{2d}(\R)$ fixed, the 
$SL_{H_1,H_2,\varepsilon}(\mathbf{f})$ converges to 
$SL_{H_1,H_2}(\mathbf{f})$. This fact combined with the characterization 
result of the convergence in $(L^{2})$ in terms of the $S$-transform, recalled 
in Proposition \ref{Prop}, allows to improve the previous statements 
concerning the intersection local times (Theorem \ref{Th11} below). In 
particular, this theorem extends the results obtained in 
\cite{NL07} to different and more general Hurst multiparameters.

\begin{proposition}\label{Prop}
Let $\left( \Phi _n\right)_{n\in\N}$ be a sequence in 
$(L^{2})$ and $\Phi\in (L^{2})$. The following two assertions are 
equivalent:
\item[{\it (i)}] $\left( \Phi _n\right)_{n\in\N}$ converges in $(L^{2})$ to 
$\Phi$;
\item[{\it (ii)}] the sequence $\left(\|\Phi _n\|\right)_{n\in\N}$ converges to $\|\Phi\|$ and, for all $\mathbf{f}\in S_{2d}(\R)$, $\left( S\Phi _n(\mathbf{f})\right)_{n\in\N}$ converges to $S\Phi(\mathbf{f})$.

Here $\|\cdot\|$ denotes the norm defined on $(L^2)$.
\end{proposition}

\begin{theorem} 
\label{Th11}
For any pair of Hurst multiparameters $H_1, H_2\in(0,1)^d$, $d\geq 1$, 
such that $d<1/{\bar H_1}+1/{\bar H_2}$, the intersection local times 
$L_{H_1,H_2}$ as well as all $L_{H_1,H_2,\varepsilon}$, 
$\varepsilon>0$, exist in $(L^{2})$, and the sequence of 
$L_{H_1,H_2,\varepsilon}$ converges in $(L^{2})$ to $L_{H_1,H_2}$ as 
$\varepsilon$ tends to zero.  
\end{theorem}

\noindent
\textbf{Proof.} According to the previous considerations, the proof amounts to 
show that $L_{H_1,H_2}, L_{H_1,H_2,\varepsilon}\in (L^2)$, for all 
$\varepsilon >0$, and that the convergence (in $\varepsilon$) of their 
$(L^2)$-norms holds. For this purpose we begin by showing that the sums
\begin{equation}
\sum_{\mathbf{m}}\sum_{\mathbf{k}}\mathbf{m}!\mathbf{k}!
\left|F_{H_1,H_2,\varepsilon, \mathbf{m},\mathbf{k}}\right|^2_{(L^2_{2d}(\R))^{\otimes (m+k)}},
\sum_{\mathbf{m}}\sum_{\mathbf{k}}\mathbf{m}!\mathbf{k}!
\left|F_{H_1,H_2,\mathbf{m},\mathbf{k}}\right|^2_{(L^2_{2d}(\R))^{\otimes (m+k)}},\label{normas}
\end{equation}
converge, where $F_{H_1,H_2,\varepsilon, \mathbf{m},\mathbf{k}}$ and 
$F_{H_1,H_2,\mathbf{m},\mathbf{k}}$ are the kernels given by 
Theorem \ref{Th9} and Proposition \ref{Proposition2}, respectively. By 
(\ref{norm}), this will prove that $L_{H_1,H_2,\varepsilon}, L_{H_1,H_2}\in(L^2)$ 
with $\|L_{H_1,H_2,\varepsilon}\|^2$ given by the first sum appearing in 
(\ref{normas}) and $\|L_{H_1,H_2}\|^2$ by the second one. 

Similar calculations done to prove Theorem \ref{Th9} yield
\begin{eqnarray*}
&&\sum_{\mathbf{m}}\sum_{\mathbf{k}}\mathbf{m}!\mathbf{k}!
\left|F_{H_1,H_2,\varepsilon, \mathbf{m},\mathbf{k}}\right|^2_{(L^2_{2d}(\R))^{\otimes (m+k)}}\\
&=&\sum_{\mathbf{m}}\sum_{\mathbf{k}}\mathbf{m}!\mathbf{k}!
\left(\frac{1}{2\pi}\right)^d\frac{(-1)^{m+3k}}{\left(\left(\frac{\mathbf{m}+\mathbf{k}}{2}\right)!\right)^2}\left(\frac12\right)^{m+k}\binom{\mathbf{m}+\mathbf{k}}{\mathbf{m}}^2\\
&&\int_0^T\!\!dt\!\!\int_0^T\!\!ds\int_0^T\!\!dt'\!\!\int_0^T\!\!ds'
\prod_{j=1}^d\left(\frac{1}{\sqrt{(\varepsilon+t^{2H_{1,j}}+s^{2H_{2,j}})(\varepsilon+{t'}^{2H_{1,j}}+{s'}^{2H_{2,j}})}}\right)^{m_j+k_j+1}\\
&&\qquad\qquad\qquad\qquad\qquad\cdot\langle M_{H_{1,j}}\I_{\left[0,t\right]},M_{H_{1,j}}\I_{\left[0,t'\right]}\rangle^{m_j}\langle M_{H_{2,j}}\I_{\left[0,s\right]},M_{H_{2,j}}\I_{\left[0,s'\right]}\rangle^{k_j},
\end{eqnarray*}
with the inner products being equal to
\begin{eqnarray*}
\langle M_{H_{1,j}}\I_{\left[0,t\right]},M_{H_{1,j}}\I_{\left[0,t'\right]}\rangle &=&
\frac{1}{2}\left(t^{2H_{1,j}}+{t'}^{2H_{1,j}}-|t-t'|^{2H_{1,j}}\right),\\
\langle M_{H_{2,j}}\I_{\left[0,s\right]},M_{H_{2,j}}\I_{\left[0,s'\right]}\rangle &=&
\frac{1}{2}\left(s^{2H_{2,j}}+{s'}^{2H_{2,j}}-|s-s'|^{2H_{2,j}}\right),
\end{eqnarray*}
for each $j=1,\ldots,d$,
\begin{eqnarray*}
&=&\left(\frac{1}{2\pi}\right)^d
\int_0^T\!\!dt\!\!\int_0^T\!\!ds\int_0^T\!\!dt'\!\!\int_0^T\!\!ds'
\prod_{j=1}^d\frac{1}{\sqrt{(\varepsilon+t^{2H_{1,j}}+s^{2H_{2,j}})(\varepsilon+{t'}^{2H_{1,j}}+{s'}^{2H_{2,j}})}}\\
&&\sum_{n=0}^\infty\frac{1}{4^nn!}\!\!\sum_{{n_1,\cdots ,n_d}\atop{n_1 +\cdots + n_d = n}}\!\frac{n!}{n_1!\ldots n_d!}\prod_{j=1}^d\frac{(2n_j)!}{n_j!}
\left(\frac{1}{(\varepsilon+t^{2H_{1,j}}+s^{2H_{2,j}})(\varepsilon+{t'}^{2H_{1,j}}+{s'}^{2H_{2,j}})}\right)^{n_j}\\
&&\cdot\left(\frac{1}{2}\left(t^{2H_{1,j}}+{t'}^{2H_{1,j}}-|t-t'|^{2H_{1,j}}+
s^{2H_{2,j}}+{s'}^{2H_{2,j}}-|s-s'|^{2H_{2,j}}\right)\right)^{2n_j}.
\end{eqnarray*}
Concerning the integrand function, observe that using the following equalities 
for the Gamma function,
\begin{eqnarray*}
&&\frac{(2n!)}{n!}=\frac{2^{2n}}{\sqrt{\pi}}\Gamma\left(n+\frac{1}{2}\right),\\
&&\Gamma\left(n+\frac{1}{2}\right)=\Gamma\left(\frac{1}{2}\right)\prod_{i=0}^{n-1}\left(\frac{1}{2}+i\right)=\sqrt{\pi}\prod_{i=0}^{n-1}\left(\frac{1}{2}+i\right),
\end{eqnarray*}
one may rewrite it as
\begin{eqnarray}
&&\prod_{j=1}^d\frac{1}{\sqrt{(\varepsilon+t^{2H_{1,j}}+s^{2H_{2,j}})(\varepsilon+{t'}^{2H_{1,j}}+{s'}^{2H_{2,j}})}}\cdot\nonumber\\
&&\cdot\sum_{n=0}^\infty\sum_{{n_1,\cdots ,n_d}\atop{n_1 +\cdots + n_d = n}}\!\prod_{j=1}^d\frac{\Gamma\left(n_j+\frac{1}{2}\right)}{\sqrt{\pi}n_j!}
\left(\frac{1}{(\varepsilon+t^{2H_{1,j}}+s^{2H_{2,j}})(\varepsilon+{t'}^{2H_{1,j}}+{s'}^{2H_{2,j}})}\right)^{n_j}\nonumber\\
&&\cdot\left(\frac{1}{2}\left(t^{2H_{1,j}}+{t'}^{2H_{1,j}}-|t-t'|^{2H_{1,j}}+
s^{2H_{2,j}}+{s'}^{2H_{2,j}}-|s-s'|^{2H_{2,j}}\right)\right)^{2n_j}\nonumber\\
&=&\prod_{j=1}^d\frac{1}{\sqrt{(\varepsilon+t^{2H_{1,j}}+s^{2H_{2,j}})(\varepsilon+{t'}^{2H_{1,j}}+{s'}^{2H_{2,j}})}}\sum_{n=0}^\infty\sum_{{n_1,\cdots ,n_d}\atop{n_1 +\cdots + n_d = n}}\prod_{j=1}^d\frac{1}{n_j!}\left(\prod_{i=0}^{n_j-1}\left(\frac{1}{2}+i\right)\right)\cdot\nonumber\\
&&\cdot\prod_{j=1}^d
\frac{\left(t^{2H_{1,j}}+{t'}^{2H_{1,j}}-|t-t'|^{2H_{1,j}}+
s^{2H_{2,j}}+{s'}^{2H_{2,j}}-|s-s'|^{2H_{2,j}}\right)^{2n_j}}{4^{n_j}(\varepsilon+t^{2H_{1,j}}+s^{2H_{2,j}})^{n_j}(\varepsilon+{t'}^{2H_{1,j}}+{s'}^{2H_{2,j}})^{n_j}}.\label{soma}
\end{eqnarray}

Hence, taking into account that for any $0\leq t,t',s,s'\leq T$ and for any 
$j=1,\ldots,d$
\[
0\leq\frac{\left(t^{2H_{1,j}}+{t'}^{2H_{1,j}}-|t-t'|^{2H_{1,j}}+
s^{2H_{2,j}}+{s'}^{2H_{2,j}}-|s-s'|^{2H_{2,j}}\right)^2}{4(\varepsilon+t^{2H_{1,j}}+s^{2H_{2,j}})(\varepsilon+{t'}^{2H_{1,j}}+{s'}^{2H_{2,j}})}<1,
\]
one recognizes that the sum in (\ref{soma}) is indeed the Taylor expansion of 
the function
\[
\prod_{j=1}^d\left(\sqrt{1-\frac{\left(t^{2H_{1,j}}+{t'}^{2H_{1,j}}-|t-t'|^{2H_{1,j}}+
s^{2H_{2,j}}+{s'}^{2H_{2,j}}-|s-s'|^{2H_{2,j}}\right)^2}{4(\varepsilon+t^{2H_{1,j}}+s^{2H_{2,j}})(\varepsilon+{t'}^{2H_{1,j}}+{s'}^{2H_{2,j}})}}\right)^{-1},
\]
and thus 
\begin{eqnarray*}
&&\sum_{\mathbf{m}}\sum_{\mathbf{k}}\mathbf{m}!\mathbf{k}!
\left|F_{H_1,H_2,\varepsilon, \mathbf{m},\mathbf{k}}\right|^2_{(L^2_{2d}(\R))^{\otimes (m+k)}}\\
&=&\left(\frac{1}{2\pi}\right)^d
\int_0^T\!\!dt\!\!\int_0^T\!\!ds\int_0^T\!\!dt'\!\!\int_0^T\!\!ds'\prod_{j=1}^d
\left((\varepsilon\!+\!t^{2H_{1,j}}\!+\!s^{2H_{2,j}})(\varepsilon\!+\!{t'}^{2H_{1,j}}\!+\!{s'}^{2H_{2,j}})\right.\\
&&\left.-\frac{1}{4}\!\left(t^{2H_{1,j}}\!+{t'}^{2H_{1,j}}\!-\!|t-t'|^{2H_{1,j}}\!+\!
s^{2H_{2,j}}\!+\!{s'}^{2H_{2,j}}\!-\!|s-s'|^{2H_{2,j}}\right)^2\right)^{-\frac{1}{2}}.
\end{eqnarray*}
Independently of the dimension $d$ and the Hurst multiparameters under 
consideration, clearly such a multiple integral is always finite.

Concerning the second sum in (\ref{normas}), first we note that the expression 
of the kernels $F_{H_1,H_2,\mathbf{m},\mathbf{k}}$ coincides with 
$F_{H_1,H_2,\varepsilon, \mathbf{m},\mathbf{k}}$ for $\varepsilon=0$. Thus, one 
may apply the previous scheme, just replacing $\varepsilon$ by zero, 
with the slight difference that in this case one has
\[
0<\frac{\left(t^{2H_{1,j}}+{t'}^{2H_{1,j}}-|t-t'|^{2H_{1,j}}+
s^{2H_{2,j}}+{s'}^{2H_{2,j}}-|s-s'|^{2H_{2,j}}\right)^2}{4(t^{2H_{1,j}}+s^{2H_{2,j}})({t'}^{2H_{1,j}}+{s'}^{2H_{2,j}})}<1,
\]
only for $0<t,t',s,s'\leq T$ such that $t\not=t'$ and $s\not=s'$. Thus, only 
for $0<t,t',s,s'\leq T$ such that $t\not=t'$ and $s\not=s'$ the sum 
corresponding to the sum in (\ref{soma}) converges. As a result, in this case 
we obtain
\begin{eqnarray*}
&&\sum_{\mathbf{m}}\sum_{\mathbf{k}}\mathbf{m}!\mathbf{k}!
\left|F_{H_1,H_2,\mathbf{m},\mathbf{k}}\right|^2_{(L^2_{2d}(\R))^{\otimes (m+k)}}\\
&=&4\left(\frac{1}{2\pi}\right)^d
\int_0^T\!\!dt\!\!\int_0^t\!\!dt'\int_0^T\!\!ds\!\!\int_0^s\!\!ds'\prod_{j=1}^d
\left((t^{2H_{1,j}}\!+\!s^{2H_{2,j}})({t'}^{2H_{1,j}}\!+\!{s'}^{2H_{2,j}})\right. \\
&&\left.-\frac{1}{4}\!\left(t^{2H_{1,j}}\!+{t'}^{2H_{1,j}}\!-\!(t-t')^{2H_{1,j}}\!+\!
s^{2H_{2,j}}\!+\!{s'}^{2H_{2,j}}\!-\!(s-s')^{2H_{2,j}}\right)^2\right)^{-\frac{1}{2}}.
\end{eqnarray*}
Due to the existence of singular points, an additional analysis is now needed 
in order to show the convergence of this multiple integral. As before, it is 
enough to consider the case $T=1$. 

As a first step we use the fact that for each $j=1,\ldots,d$ fixed one has
\begin{eqnarray}
&&(t^{2H_{1,j}}\!+\!s^{2H_{2,j}})({t'}^{2H_{1,j}}\!+\!{s'}^{2H_{2,j}})\label{original}\\
&&-\frac{1}{4}\!\left(t^{2H_{1,j}}\!+{t'}^{2H_{1,j}}\!-\!|t-t'|^{2H_{1,j}}\!+\!
s^{2H_{2,j}}\!+\!{s'}^{2H_{2,j}}\!-\!|s-s'|^{2H_{2,j}}\right)^2\nonumber\\
&\geq& t^{2H_{1,j}}{t'}^{2H_{1,j}}-\frac{1}{4}\left(t^{2H_{1,j}}+{t'}^{2H_{1,j}}-|t-t'|^{2H_{1,j}}\right)^2\label{function1}\\
&&+s^{2H_{2,j}}{s'}^{2H_{2,j}}-\frac{1}{4}\left(s^{2H_{2,j}}+{s'}^{2H_{2,j}}-|s-s'|^{2H_{2,j}}\right)^2,\label{function2}
\end{eqnarray}
with the advantage that, in contrast to (\ref{original}), (\ref{function1}) 
as well as (\ref{function2}) only depend of a unique Hurst parameter. 
Moreover, (\ref{function1}) and (\ref{function2}) are both of the type
\[
\varphi_H(u,v):=u^{2H}{v}^{2H}-\frac{1}{4}\left(u^{2H}+v^{2H}-|u-v|^{2H}\right)^2,
\]
which, as a function of $u$ and $v$, is an homogeneous function of order $4H$. 
Therefore, for every $0<v<u<1$ one has
\begin{eqnarray*}
&&u^{2H}{v}^{2H}-\frac{1}{4}\left(u^{2H}+v^{2H}-(u-v)^{2H}\right)^2\\
&=&u^{4H}\left[\left(\frac{v}{u}\right)^{2H}-\frac{1}{4}\left(1+\left(\frac{v}{u}\right)^{2H}-\left(1-\frac{v}{u}\right)^{2H}\right)^2\right],
\end{eqnarray*}
where, for $v/u\in\left(0,1\right)$ fixed, the expression between the square 
brackets is a decreasing function of $H\in\left(0,1\right)$. 

As a consequence,
\begin{eqnarray}
&&\sum_{\mathbf{m}}\sum_{\mathbf{k}}\mathbf{m}!\mathbf{k}!
\left|F_{H_1,H_2,\mathbf{m},\mathbf{k}}\right|^2_{(L^2_{2d}(\R))^{\otimes (m+k)}}\nonumber\\
&\leq&4\left(\frac{1}{2\pi}\right)^d
\int_0^1\!\!dt\!\!\int_0^t\!\!dt'\int_0^1\!\!ds\!\!\int_0^s\!\!ds'
\left(t^{2\bar H_1}{t'}^{2\bar H_1}-\frac{1}{4}\left(t^{2\bar H_1}+{t'}^{2\bar H_1}-(t-t')^{2\bar H_1}\right)^2\right.\nonumber\\
&+&\left.s^{2\bar H_2}{s'}^{2\bar H_2}-\frac{1}{4}\left(s^{2\bar H_2}+{s'}^{2\bar H_2}-(s-s')^{2\bar H_2}\right)^2\right)^{-\frac{d}{2}}.\label{int}
\end{eqnarray}

Now the proof follows closely the one in \cite[Proof of Lemma 4]{NL07}, based 
on the fact that
\[
\lambda^{-\frac{d}{2}}=\frac{1}{\Gamma(\frac{d}{2})}\int_0^{+\infty}dz\,e^{-\lambda z}z^{\frac{d}{2}-1},
\]
which allows to rewrite the multiple integral in (\ref{int}) as
\begin{equation}
\frac{1}{\Gamma(\frac{d}{2})}\int_0^{+\infty} dz\,z^{\frac{d}{2}-1}
\left(\int_0^1dt\int_0^tds\, e^{-z\varphi_{\bar H_1}(t,s)}\right)
\left(\int_0^1dt\int_0^tds\, e^{-z\varphi_{\bar H_2}(t,s)}\right).\label{Nov}
\end{equation}
Since
\[
\forall\,z\in\left[0,1\right],\ 
\int_0^1dt\int_0^tds\, e^{-z\varphi_{\bar H_i}(t,s)}<+\infty,\quad i=1,2,
\]
the convergence of the integral (\ref{Nov}) then will follow from the 
convergence of the integral
\[
\int_1^{+\infty} dz\,z^{\frac{d}{2}-1}
\left(\int_0^1dt\int_0^tds\, e^{-z\varphi_{\bar H_1}(t,s)}\right)
\left(\int_0^1dt\int_0^tds\, e^{-z\varphi_{\bar H_2}(t,s)}\right).
\]

As in \cite[Proof of Lemma 4]{NL07}, the homogeneity property of 
$\varphi_{\bar H_1}$ and $\varphi_{\bar H_2}$ yields
\begin{eqnarray*}
&&\int_1^{+\infty}\!dz\,z^{\frac{d}{2}-1}
\left(\int_0^1dt\int_0^tds\, e^{-z\varphi_{\bar H_1}(t,s)}\right)
\left(\int_0^1dt\int_0^tds\, e^{-z\varphi_{\bar H_2}(t,s)}\right)\\
&=&\!\!\!\int_1^{+\infty} dz\,z^{\frac{d}{2}-1-\frac{1}{2\bar H_1}-\frac{1}{2\bar H_2}}
\left(\int_0^{z^{\frac{1}{4\bar H_1}}}\!\!\!dx\int_0^xdy\, e^{-\varphi_{\bar H_1}(x,y)}\right)
\left(\int_0^{z^{\frac{1}{4\bar H_2}}}\!\!\!dx\int_0^xdy\, e^{-\varphi_{\bar H_2}(x,y)}\right),
\end{eqnarray*}
where a double change of coordinates leads for each $i=1,2$ to
\begin{eqnarray*}
&&\int_0^{z^{\frac{1}{4\bar H_i}}}\!\!\!dx\int_0^xdy\, e^{-\varphi_{\bar H_i}(x,y)}\\
&\leq&\frac{1}{4\bar H_i}\int_0^{\pi/4}d\theta\,(\varphi_{\bar H_i}(\cos\theta,\sin\theta))^{-\frac{1}{2\bar H_i}}\gamma\left(\frac{1}{2\bar H_i},2^{2\bar H_i}z\varphi_{\bar H_i}(\cos\theta,\sin\theta)\right).
\end{eqnarray*}
Here $\gamma$ is the lower incomplete gamma function, that is,
\[
\gamma(\alpha,x):=\int_0^xdy\,e^{-y}y^{\alpha-1},\quad \alpha>0,
\]
which, as shown in \cite[Lemma 2]{NL07}, is bounded by
\[
\gamma(\alpha,x)\leq K(\alpha)x^\epsilon,\quad 
K(\alpha):=\max\left\{\frac{1}{\alpha},\Gamma(\alpha)\right\},
\]
for all $x>0$ and for every $0<\epsilon<\alpha$. 

Hence, for all $0<\epsilon<\frac{1}{2\max\{\bar H_1,\bar H_2\}}$ one 
finally obtains
\begin{eqnarray}
&&\int_1^{+\infty}\!dz\,z^{\frac{d}{2}-1}
\left(\int_0^1dt\int_0^tds\, e^{-z\varphi_{\bar H_1}(t,s)}\right)
\left(\int_0^1dt\int_0^tds\, e^{-z\varphi_{\bar H_2}(t,s)}\right)\nonumber\\
&\leq&\frac{2^{2\epsilon(\bar H_1+\bar H_2)}}{16\bar H_1\bar H_2}K\left(\frac{1}{2\bar H_1}\right)K\left(\frac{1}{2\bar H_2}\right)\int_1^{+\infty} dz\,z^{\frac{d}{2}-1-\frac{1}{2\bar H_1}-\frac{1}{2\bar H_2}+2\epsilon}\cdot\label{Err}\\
&&\cdot\left(\int_0^{\pi/4}d\theta\,(\varphi_{\bar H_1}(\cos\theta,\sin\theta))^{\epsilon-\frac{1}{2\bar H_1}}\right)\left(\int_0^{\pi/4}d\theta\,(\varphi_{\bar H_2}(\cos\theta,\sin\theta))^{\epsilon-\frac{1}{2\bar H_2}}\right).\nonumber
\end{eqnarray}
Concerning the integral in $z$, clearly it converges provided 
$\epsilon < \frac{\bar H_1+\bar H_2 -d\bar H_1\bar H_2}{4\bar H_1\bar H_2}$, being $\frac{\bar H_1+\bar H_2 -d\bar H_1\bar H_2}{4\bar H_1\bar H_2}$ 
always a positive number, because $d<1/{\bar H_1}+1/{\bar H_2}$. These facts 
combined mean that in (\ref{Err}) one shall fix a
\[
0<\epsilon<\min\left\{\frac{1}{2\max\{\bar H_1,\bar H_2\}},\frac{\bar H_1+\bar H_2 -d\bar H_1\bar H_2}{4\bar H_1\bar H_2}\right\}.
\] 
For such a $\epsilon$ fixed, also both integrals 
in $\theta$ converge, cf.~\cite[Proof of Lemma 4]{NL07}, and thus (\ref{Err}) 
converges.

In this way we have shown that, on the one hand, $L_{H_1,H_2}\in (L^2)$, and, 
on the other hand, that one may apply a Lebesgue dominated convergence 
argument to infer the convergence in $\varepsilon$ of
$\|L_{H_1,H_2,\varepsilon}\|^2$ to $\|L_{H_1,H_2}\|^2$.\hfill$\blacksquare \medskip$

\begin{remark}
Under the conditions of Theorem \ref{Th11}, one has
\[
\min\left\{\frac{1}{2\max\{\bar H_1,\bar H_2\}},\frac{\bar H_1+\bar H_2 -d\bar H_1\bar H_2}{4\bar H_1\bar H_2}\right\}=\frac{\bar H_1+\bar H_2 -d\bar H_1\bar H_2}{4\bar H_1\bar H_2}
\]
whenever 
$d\geq\frac{1}{\min\{\bar H_1,\bar H_2\}}-\frac{1}{\max\{\bar H_1,\bar H_2\}}$, while for 
$d<\frac{1}{\min\{\bar H_1,\bar H_2\}}-\frac{1}{\max\{\bar H_1,\bar H_2\}}$,
\[
\min\left\{\frac{1}{2\max\{\bar H_1,\bar H_2\}},\frac{\bar H_1+\bar H_2 -d\bar H_1\bar H_2}{4\bar H_1\bar H_2}\right\}=\frac{1}{2\max\{\bar H_1,\bar H_2\}}.
\]
\end{remark}

\section*{Appendix}

\begin{lemma}
\label{Lemma1} Let $0<H_1<H_2<1$ be given. The integral
\[
\int_0^1dt\int_0^1ds\,\frac{(t^2+s^2)^N}{(t^{2H_1}+s^{2H_2})^{N+\frac{d}{2}}}
\] 
is finite if and only if $2H_2\left(N+\frac{d}{2}\right)<1+2N+\frac{H_2}{H_1}$.
\end{lemma} 

\noindent
\textbf{Proof.} Through the change of variables $t=u^{\frac{H_2}{H_1}}$ one 
obtains
\begin{eqnarray}
&&\int_0^1dt\int_0^1ds\,\frac{(t^2+s^2)^N}{(t^{2H_1}+s^{2H_2})^{N+\frac{d}{2}}}\nonumber\\
&=&\frac{H_2}{H_1}\sum_{n=0}^N\binom{N}{n}\int_0^1du\int_0^1ds\,
\frac{u^{2n\frac{H_2}{H_1}}s^{2(N-n)}}{(u^{2H_2}+s^{2H_2})^{N+\frac{d}{2}}}
u^{\frac{H_2}{H_1}-1},\label{integral}
\end{eqnarray}
where each double integral appearing in (\ref{integral}) is finite if and only 
if the integrand function is integrable on the unit ball $B_1(0)\subset\R^2$. 
For each $n=0,1,...,N$ a polar change of coordinates yields
\begin{eqnarray*}
&&\int_{B_1(0)}duds\,\frac{u^{2n\frac{H_2}{H_1}}s^{2(N-n)}}{(u^{2H_2}+s^{2H_2})^{N+\frac{d}{2}}}u^{\frac{H_2}{H_1}-1}\\
&=&\int_0^{2\pi}d\theta\,\frac{\cos^{(2n+1)\frac{H_2}{H_1}-1}\theta\cdot\sin^{2(N-n)}\theta}{\left(\cos^{2H_2}\theta+\sin^{2H_2}\theta\right)^{N+\frac{d}{2}}}\int_0^1dr\,\frac{1}{r^{2H_2(N+\frac{d}{2})-(2n+1)\frac{H_2}{H_1}-2(N-n)}},
\end{eqnarray*}
which is finite if and only if the integral in $r$ converges, that is, if and 
only if $2H_2(N+\frac{d}{2})-(2n+1)\frac{H_2}{H_1}-2(N-n)<1$. This shows that a necessary and sufficient condition for the convergence of the sum in 
(\ref{integral}) is given by $2H_2(N+\frac{d}{2})-2N-\frac{H_2}{H_1}<1$.\hfill$\blacksquare \medskip$

\begin{lemma}
\label{Lemma2} Given $H_i=(H_{i,1},...,H_{i,d})\in\left(0,1\right)^d$, $i=1,2$, 
assume that $\bar H_1<\bar H_2$. Then
\[
\int_0^1dt\int_0^1ds\,\prod_{j=1}^d\frac{1}{\sqrt{t^{2H_{1,j}}+s^{2H_{2,j}}}}\left(\frac{t^2+s^2}{t^{2\bar H_1}+s^{2\bar H_2}}\right)^N<\infty
\]
whenever $2\bar H_2\left(N+\frac{d}{2}\right)<1+2N+\frac{\bar H_2}{\bar H_1}$.
\end{lemma}

\noindent
\textbf{Proof.} Since in $\left[0,1\right]^2$ the following inequality holds
\[
\prod_{j=1}^d\frac{1}{\sqrt{t^{2H_{1,j}}+s^{2H_{2,j}}}}\left(\frac{t^2+s^2}{t^{2\bar H_1}+s^{2\bar H_2}}\right)^N\leq 
\frac{\left(t^2+s^2\right)^N}{\left(t^{2\bar H_1}+s^{2\bar H_2}\right)^{N+\frac{d}{2}}},
\]
the proof reduces to an application of the previous lemma.\hfill$\blacksquare \medskip$

\subsection*{Acknowledgments}

We truly thank F.~P.~da Costa for the helpful discussions. Financial support 
of projects PTDC/MAT/67965/2006, PTDC/MAT/100983/2008 and FCT, POCTI-219, 
ISFL-1-209 is gratefully acknowledged.


\end{document}